\documentclass[twocolumn]{amsart}

\usepackage{amssymb,amsmath,amscd,graphicx,xfrac}

\usepackage{tikz}
\tikzset{every picture/.style={line width=0.75pt}}
\usetikzlibrary{shapes.geometric}
\usetikzlibrary{quotes,angles,positioning}

\DeclareMathAlphabet{\mymathbb}{U}{bbold}{m}{n}

\newcommand{\XX}{\mathbb{X}}
\newcommand{\YY}{\mathbb{Y}}
\newcommand{\NN}{\mathbb{N}}
\newcommand{\RR}{\mathbb{R}}
\newcommand{\ZZ}{\ts\mathbb{Z}}
\newcommand{\QQ}{\mathbb{Q}}
\newcommand{\CC}{\mathbb{C}}
\newcommand{\cA}{\mathcal{A}}
\newcommand{\cL}{\mathcal{L}}
\newcommand{\cT}{\mathcal{T}}

\newcommand{\vL}{\varLambda}

\newcommand{\ts}{\hspace{0.5pt}}
\newcommand{\nts}{\hspace{-0.5pt}}

\theoremstyle{definition}

\newcommand{\ee}{\ts\mathrm{e}}
\newcommand{\dd}{\,\mathrm{d}\ts}

\newcommand{\dens}{\mathrm{dens}}
\newcommand{\vol}{\mathrm{vol}}
\newcommand{\ii}{\ts\mathrm{i}}

\renewcommand{\epsilon}{\varepsilon}

\newcommand{\oplam}{\mbox{\Large $\curlywedge$}}

\newcommand{\defeq}{\mathrel{\mathop:}=}

\newcommand{\myfrac}[2]{\frac{\raisebox{-2pt}{$#1$}}
	{\raisebox{0.5pt}{$#2$}}}

\begin{document}	
\twocolumn[
\begin{LARGE}
	\centerline{On the Fibonacci tiling and its
		modern ramifications}\vspace{3ex}
\end{LARGE}
\centerline{\large Michael Baake, Franz G\"{a}hler and Jan
  Maz\'{a}\v{c}}\vspace{2ex}
\begin{footnotesize}
  \centerline{\textit{Fakult\"at f\"ur Mathematik,
      Universit\"at Bielefeld, Postfach 100131, 33501 Bielefeld,
      Germany\vspace{1mm}}}
\end{footnotesize}\vspace{4ex}
\begin{small}
\hrule\vspace{2ex}
\begin{minipage}{\textwidth}
  \textbf{Abstract}\vspace{2ex}\\
  In the last 30 years, the mathematical theory of aperiodic order has
  developed enormously. Many new tilings and properties have been
  discovered, few of which are covered or anticipated by the early
  papers and books. Here, we start from the well-known Fibonacci chain
  to explain some of them, with pointers to various generalisations as
  well as to higher-dimensional phenomena and results. This should
  give some entry points to the modern literature on the subject.
  \vspace{2ex}\\

	\end{minipage}\vspace{2ex}
	\hrule
\end{small}\vspace{6ex}
]

\section{Introduction}\label{sec:intro}

Let us begin with a rough sketch of the different perspectives on
the field of aperiodic order, and how it developed.  In mathematics,
the origins of the field of aperiodic order are two-fold. On the one
hand, the connection between non-decidability questions and the
existence of aperiodic tile sets was instrumental to the investigation
of aperiodic tilings. On the other hand, the theory of almost-periodic
functions due to Harald Bohr \cite{Bohr} showed the existence of
long-range order beyond periodicity, though this was little known or
appreciated in the physical sciences.

In crystallography and physics, the detailed study of structural
disorder and incommensurate phenomena slowly paved the ground for
going beyond ordinary crystals, which clearly showed the need for an
extension of classic solid state physics.  Shechtman's discovery
\cite{SBGC84} of icosahedral quasicrystals in 1982 then started a
rapid development in many directions, both mathematical and physical.

Initially, the effort in mathematics and physics was largely
synchronous, clearly driven by many open questions and the need for
new tools to answer them. After a while, crystallography and physics
became largely satisfied with the new toolbox, even though the
connection between the theoretical and the more applied branches still
seemed somewhat speculative in places.

In particular, despite the success of tilings of Penrose type in the
description of quasiperiodic long-range order, no clear connection
between aperiodic tile sets and real-world quasicrystals could be
established (which still is the case today, though some new phenomena
connected with monotiles might change this now). At the same time,
mathematicians wanted to explore these structures and their
possibilities without any real-world constraints.  Therefore, the
mathematical and physical research directions gradually drifted apart
and followed their own goals, as is often the case after a decade (or
so) of joint effort.

Since the mid-1990s, the mathematical theory of aperiodic order really
took off at an amazing pace, and rather little of the outcome was
noticed in the physical sciences. Likewise, only some mathematicians
kept an eye on new results in the (experimental) quasicrystal
world. Each side developed new methods and produced results relevant
to the other one, but the impact on one another, unfortunately, was
relatively small. This also concerned the connection between aperiodic
tile sets and quasicrystals.

Some progress then came, quite unexpectedly, via a~decorated hexagonal
tile, originally due to Joan Taylor and then further analysed by
Socolar and Taylor \cite{ST}, which was a functional monotile (when
one also admits its reflected version) \emph{and} a mathematical
quasicrystal; see \cite{TAO,ML}. It had no purely geometric
realisation with a disk-like tile though, but needed nearest and
next-to-nearest neighbour information to encode perfect local (or
matching) rules.

This situation recently changed with the discovery of the \emph{Hat}
family of monotiles \cite{Hat}. Each of them enforces aperiodicity by
a purely geometric face-to-face condition, yet also with the need to
admit the reflected version. Soon after, the same author team
constructed a completely chiral analogue, now known as the
\emph{Spectre} \cite{Spectre}. Both define aperiodic tilings with
disk-like tiles in a purely geometric way, but they do not specify a
unique LI class (see below for more on this notion) of tilings.

It was one declared goal of ICQ15 to bring the mathematical and
physical sides together again. In mathematics, to which this little
survey concentrates, two major directions were identified, namely
\emph{topological structures} and \emph{invariants} on the one side
\cite{Johannes} and \emph{diffraction} and \emph{spectral theory} on
the other. Our main problem now is to summarise three decades of
mathematical development in an introductory way that does not assume
knowledge of all recent methods, which seems impossible. As a
compromise, we attempt to start from the best-studied one-dimensional
example, the Fibonacci chain, and describe as many aspects as possible
on the basis of it. The concepts explained below also apply to a
higher-dimensional setting, possibly with some small adjustments, and
then cover all the famous examples such as the Penrose, the
Ammann--Beenker, or the shield and the square-triangle tilings, and
many more; see \cite{TAO,tilings} for details.

Amazingly, our mainly one-dimensional approach connects to quite a few
modern and recent results, though we will often be sketchy and refer
to the relevant literature for details. Most terms are defined and
illustrated more extensively in the monograph \cite{TAO}, while rather
little is coverend in the early literature. Some other technical terms
can easily be looked up in the {\sc WikipediA}, which has turned into
a decent source for mathematical definitions and explanations. We thus
assume that the reader will use these two sources, augmented by the
references provided here. We hope that this will provide a basis to
allow an entry point or even some stepping stones to the present state
of the art in aperiodic order.

Clearly, the Fibonacci example is not always sufficient. For instance,
the above-mentioned Hat and Spectre tilings are truly two-dimensional
affairs. Nevertheless, even some aspects of them can be better
understood with the tools and methods explainable for the Fibonacci
chain. Using them, we show that they are examples of structures with
pure point (or Bragg) diffraction and thus bring the original strands
together --- in an unexpected way. \smallskip

This paper is organised as follows. We first set the scene, in
Section~\ref{sec:scene}, by recalling the basic steps to generate the
Fibonacci chain and tiling, which is then followed by the projection
description, where we show how the embedding is made from intrinsic
data (Section~\ref{sec:embedding}). Some variations and complications
are discussed in Section~\ref{sec:compli}, before we sketch the
equidistribution properties of the system in
Section~\ref{sec:ergodic}. This is often tacitly assumed, but far from
trivial, and it is the basis for practically all ergodic arguments
used in averaging over the (infinite) system.  We then dive into
various aspects of the pair correlations (Section~\ref{sec:pair}),
which have a nice dual interpretation --- namely from the embedding
picture and via an exact renormalisation scheme.

This is followed by a summary of the possible shape changes in
Section~\ref{sec:shape}, and how it can be understood in the
projection approach.  With this, we are prepared for a discussion of
the diffraction properties of the Fibonacci chain and its variants
(Section~\ref{sec:diffract}), which mimics the situation one has to
face in the recently discovered monotile tilings.  Finally, we sketch
the dynamical systems approach in Section~\ref{sec:dynamical}, which
is instrumental in much if not most of the recent progress in the
mathematics of aperiodic order. While we go along, we mention various
extensions and higher-dimensional analogues, with references to recent
or neglected papers.

\section{Setting the scene}\label{sec:scene}

The binary Fibonacci substitution, say on the alphabet $\cA=\{a,b\}$,
is arguably the most frequently studied one. It comes in two versions,
\begin{equation}\label{eq:Fibo-sub-def}
  \rho^{}_1  \, = \, (ab, a)
  \quad \mbox{and} \quad
  \rho^{}_2  \, = \, (ba, a) \ts ,
\end{equation}
where we write the substitution by simply listing the images of $a$
and $b$. The iteration of the legal seed $a|a$ under $\rho^{}_{1}$
gives
\begin{equation}
\begin{split}  
  a|a & \rightarrow ab|ab \rightarrow aba|aba \rightarrow abaab|abaab
  \\ & \rightarrow \cdots \rightarrow w \rightarrow w' = \rho^{}_1(w)
  \rightarrow w \rightarrow \cdots
\end{split}
\end{equation}
leading to a $2$-cycle of bi-infinite words. We write them as
$w = \ldots w^{}_{-2} w^{}_{-1} | w^{}_{0} w^{}_{1} \ldots $, which
explains the role of the marker. The words $w$ and $w'$ differ only at
the two positions immediately left of $|$, which read either $ab$ (in
$w$, say) or $ba$ (in $w'=\rho^{}_{1}(w)$).
	
Let us now take the orbit under the shift action of $\ZZ$. Here, we
define (powers of) the shift $S$ as 
\begin{equation}\label{eq:shift}
  (S^k w )^{}_n \, \defeq \,
  w^{}_{n+k} \qquad \mbox{with} \ k,n\in\ZZ \ts , 
\end{equation}
and then take the closure in the product topology (see below for
more). This gives the discrete or \emph{symbolic hull} of the
substitution, $\XX_0$. Here, it does not matter whether one starts
with $w$ or with $w'$, which are \emph{locally indistinguishable}
(LI), meaning that every subword of $w$ occurs in $w'$ and vice versa.
In fact, any two elements of $\XX_0$ are LI, and $\XX_0$ is the LI
class defined by $\rho^{}_{1}$.  Doing the analogous exercise with
$\rho^{}_2$ instead gives a different $2$-cycle (with the reflected
versions of $w$ and $w'$), but the same LI class, $\XX_0$.

At this point, we recall that $w$ and $w'$ are equal to the right
of~$|$ (the marker for the origin), but differ to the left. Such
a~singular pair is called \emph{proximal} (in fact, asymptotic), and
its existence in one LI class immediately implies the non-periodicity
of $w$, hence also of all elements of $\XX_0$, and thus
aperiodicity. Here, a bi-infinite sequence is called \emph{aperiodic}
when no element of its hull has any non-trivial period, see
\cite[Secs.~ 3.1 and 4.2]{TAO} for a more detailed discussion of why
this is important to be distinguished from mere nonperiodicity.
	
Generally, the \emph{hull} of a sequence is the closure of its shift
orbit in the product topology,
\begin{equation}
  \XX \, = \,  \XX (w) \, = \,
  \overline{ \{S^k w : k \in \ZZ \} } \ts .
\end{equation}
Here, the term \emph{product topology} refers to the natural topology
of the sequence space $\cA^{\ZZ}$ as induced by the discrete topology
on the alphabet $\cA$. Two sequences are close in the product (or
local) topology when they agree on a large region around the origin.
A hull is called \emph{minimal}, when it consists of a single LI
class, as in the case of our Fibonacci example.  Minimality of
$\XX(w)$ is equivalent to $w$ being \emph{repetitive}, which means
that every subword of $w$ occurs repeatedly in $w$, with bounded gaps.
	
The Fibonacci sequence is also \emph{Sturmian}, which refers to
repetitive bi-infinite words that have $n+1$ distinct subwords of
length $n$ for every $n\in\NN$.  They are binary aperiodic sequences
of minimal complexity and, thus, in this sense, the simplest aperiodic
sequences to study.
	
With Eq.~\eqref{eq:shift}, we have an action of the group $\ZZ$, via
the shift $S$ on $\XX_0$, which is continuous in the product topology,
and the pair $(\XX_0,\ZZ)$ is then a \emph{topological dynamical
  system}. It is worth noting that $\XX_0$ is inversion symmetric, in
line with our previous statement that $\rho^{}_1$ and $\rho^{}_2$ both
define $\XX_0$. This need not be the case, as one can see for the
substitution $(aab,ba)$ where we get an enantiomorphic (or mirror)
pair, the reflected version being generated by $(baa, ab)$.  We shall
return to this example later.
	
\begin{figure}
  \includegraphics[width=0.45\textwidth]{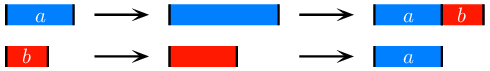}
  \caption{\small The geometric inflation rule for the
    Fibonacci chain, with natural tile lengths. \label{fig:Fibo}}
\end{figure}

A powerful quantity attached to a substitution rule $\rho$ is its
\emph{substitution matrix} $M_\rho = \bigl(m_{ij}\bigr)$, where
$m_{ij}$ is the number of letters of type $i$ in a superword of type
$j$, i.e., in $\rho(j)$ (after fixing a numbering of the
alphabet). Both Fibonacci substitutions, $\rho^{}_1$ and $\rho^{}_2$,
share the matrix
\begin{equation}\label{eq:matM}
  M \, = \, \begin{pmatrix} 1 & 1 \\ 1 & 0 \end{pmatrix},
\end{equation}
with eigenvalues $\lambda_{\pm} = \tfrac{1\pm\sqrt{5}}{2}$. The
leading one is the \emph{Perron--Frobenius} (PF) eigenvalue (the golden
ratio $\tau$ in this case). The corresponding left and right
eigenvectors are
\begin{equation}
  \langle u | \, = \, (\tau,1) \quad  \text{and} \quad
   |v\rangle \, = \, ( \tau^{-1}, \tau^{-2} )^{\top},
\end{equation} 
which are suitably normalised according to their meaning. The entries
of $|v\rangle$ encode the relative frequencies of the letters $a$ and
$b$ in the (bi-)infinite Fibonacci words, while those of $\langle u|$
give the natural tile lengths in the induced (geometric)
\emph{inflation rule}, where stretched versions of the intervals
(tiles) are subdivided according to the substitution rule as
illustrated in Figure~\ref{fig:Fibo}. In this version, we have
$\langle u \ts | \ts v \rangle = 3-\tau =\sqrt{5}/\tau$, which is the
average length of an interval in our Fibonacci tiling.

This way, one creates tilings $\cT$ of the real line that give rise to
a \emph{tiling space} or \emph{hull}
\begin{equation}
  \YY \, = \,  \overline{\{t+\mathcal{T}  :  t\in\RR  \}} \ts ,
\end{equation}
where the closure is taken in the \emph{local topology}. Here, two
tilings are $\epsilon$-close if they agree on the interval
$\bigl[-\tfrac{1}{\epsilon},\tfrac{1}{\epsilon}\bigr]$, possibly after
a small translation of one of them by at most $\epsilon$. This gives
the topological dynamical system $(\YY,\RR)$, which is the continuous
counterpart of the shift space $(\XX_0,\ZZ)$ from above.
	
Any such tiling can be turned into a \emph{Delone set} (a point set
that is both uniformly discrete and relatively dense, see \cite{TAO}
for details) in many ways, perhaps the most common one emerging from
taking the left endpoints of the intervals. They can be coloured if
tiles of equal length have to be distinguished, as happens in all
examples with degeneracies in the left PF eigenvector.

The tiling and the Delone set of control points (possibly coloured, if
one needs to distinguish control points of different tiles with the
same geometry) are locally equivalent in the sense that a strictly
local rule exists to turn them into one another. They are thus
\emph{mutually locally derivable} (MLD). Such rules commute with the
translation action and constitute an important subclass of
\emph{topological conjugacies}, which are homeomorphisms of the hull
that commute with the translation action but need not stem from a
local rule. This distinction is characteristic for aperiodic tiling
spaces, compare \cite{BGM-Fib}, and will become particularly relevant in
the context of shape changes.

\section{Embedding method}\label{sec:embedding}	

\begin{figure}
\includegraphics[width=0.45\textwidth]{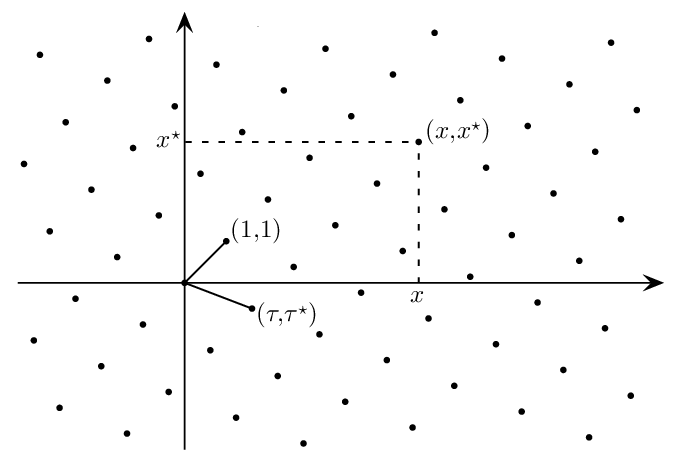}
\caption{\small The Minkowski embedding of $\ZZ[\tau]$ as a planar
  lattice, where the projections are illustrated for the lattice point
  $(x,x^{\star})$ with $x=4+\tau$. The horizontal (vertical) axis
  represents the physical (internal) space. The grey strip contains
  all lattice points that are projected to the horizontal line for our
  guiding Fibonacci example.}
\label{fig:ztau}
\end{figure}

Returning to our guiding example, now in the form of a Fibonacci
tiling (or point set), we look at the important set of translations
that shift a tile (or a patch of tiles) to another occurrence within
the same tiling, the so-called \emph{return vectors}. Due to the
inflation structure, we can do this for patches of arbitrary size in
one step. Now, also taking all integer linear combinations of return
vectors complete them into the \emph{return module}, which is
\begin{equation}
  \ZZ[\tau] \, = \, \{m+n\tau  :  m,n\in\ZZ \}
  \, = \, \ZZ \oplus \tau\ZZ
\end{equation}
for our Fibonacci example. This is a $\ZZ$-module of rank $2$ that is
a dense subset of $\RR$. It can be seen as a~projection of a lattice
in $\RR^2$ in many ways. A~particularly natural one emerges from the
\emph{Minkowski embedding} as follows. Recall that algebraic
conjugation in the quadratic field $\QQ(\sqrt{5})$ is given by sending
$\sqrt{5}$ to $-\sqrt{5}$. Denoting the corresponding mapping by
$(.)^{\star}$, we can consider
\begin{equation}
  \begin{split}
  \cL \, & = \, \{(x,x^{\star}) : x\in \ZZ[\tau] \} \\
  & = \, \ZZ (1,1) \oplus \ZZ(\tau, 1-\tau) \ts ,
\end{split}
\end{equation}
which is the lattice shown in Figure~\ref{fig:ztau}.

Note that this lattice is neither a square nor a~rectangular lattice,
but by construction satisfies that its projection to the horizontal
axis (the first coordinate) is $\ZZ[\tau]$. If one insists, a
rescaling of the second coordinate can be used to turn $\cL$ into
a~square lattice, compare \cite[Rem.~3.4]{TAO}, but this brings no
particular advantage. This just reflects the fact that the scale of
internal space relative to physical space is arbitrary, and not a
physically meaningful quantity.

The significant aspect of the return module and its embedding is the
alternative description of the Fibonacci point set as a lattice
projection set. If $\vL$ and $\vL'$ are the Delone sets of the
Fibonacci $2$-cycle for $\rho^{}_{1}$ from Eq.~\ref{eq:Fibo-sub-def},
one finds
\begin{equation}
  \vL^{(\prime)} \, = \, \oplam \bigl(W^{(\prime)}\bigr) \, \defeq \,
  \{x\in\ZZ[\tau] : x^{\star} \in W^{(\prime)} \}
\end{equation} 
where $W = (-1, \tau-1]$ and $W' = [-1,\tau-1)$ are then called the
\emph{windows} for the projection point sets $\vL$ and $\vL'$,
respectively. When the window is bounded and has a non-empty interior,
the set $\oplam (W)$ is called a \emph{model set}, and it is called
\emph{regular} when the boundary of $W$ is of zero measure; see
\cite{TAO} for details.

In our guiding example, the difference between $\vL$ and $\vL'$
consists in only two points, namely $-1$ (which belongs to $\vL'$ but
not to $\vL$) and $-\tau$ (which lies in $\vL$ but not in $\vL'$). The
union $\vL \cup \vL'$ is coded by the closed interval
$[-1,\tau-1] = W \cup W'$ as window, which is the closure of either
contributing window. Here, $\vL$ and $\vL'$ constitute an
\emph{asymptotic pair} with a difference only occuring in a bounded
region (in higher dimensions, one can have differences along
subspaces, as one knows from the Penrose worms). The existence of
asymptotic pairs is an essential feature for aperiodic, repetitive
tilings. In fact, they are the reason why Bohr's theory of almost
periodic functions needs to be extended to cover these examples, as
explained in some detail in \cite{LSS}.

There are several equivalent ways to view and interpret this
projection approach in more generality. Here, a Euclidean
\emph{cut-and-project scheme} (CPS) is a triple $(\RR^n,\RR^m,\cL)$
with a lattice $\cL \subset \RR^{n+m}$ and two natural projections
$\pi:\RR^{n+m} \rightarrow \RR^{n}$ and
$\pi^{}_{\text{int}}:\RR^{n+m} \rightarrow \RR^{m}$ satisfying that
$\pi\bigl|_{\cL}$ is injective and $\pi^{}_{\text{int}}(\cL)$ is dense
in $\RR^{m}$. Since the projection $\pi$ restricted to the lattice
$\cL$ provides a bijection between the lattice $\cL$ and
$L= \pi(\cL)$, one defines the \emph{star map}
$\star: \RR^n \rightarrow \RR^m$ as
\begin{equation}
  x \mapsto x^{\star} \defeq \pi^{}_{\text{int}} \bigl(\bigl(
  \pi\bigl|_{\cL} \bigr)^{-1}(x) \bigr).
\end{equation}
This gives a rather natural connection between the \emph{physical
  space} $\RR^n$ and the \emph{internal space} $\RR^m$. Once a CPS has
been fixed, one can choose a (sufficiently nice) subset
$W \subset \RR^m$ and consider the model set $\oplam (W)$ in complete
analogy to above.

In our guiding Fibonacci example, we have the following Euclidean CPS.
\begin{equation}\label{eq:CPS}
  \renewcommand{\arraystretch}{1.2}\begin{array}{c@{}c@{}c@{}c@{}c@{}l}
       \RR & \;\, \xleftarrow{\;\;\; \pi \;\;\; }
         & \; \RR \nts\nts \times \nts\nts \RR \! & 
         \xrightarrow{\;\: \pi^{}_{\text{int}} \;\: } & \RR & \\
          \cup & & \; \cup & &  \cup  & \hspace*{-3.5ex} 
	 \raisebox{1pt}{\text{\scriptsize dense}} \\
	 \pi (\cL) & \;\, \xleftarrow{\;\ts 1-1 \;\ts } & \; \cL & 
	 \xrightarrow{ \qquad } &\pi^{}_{\text{int}} (\cL) & \\
	 \| & & & & \| & \\
         \, L = \ZZ[\tau] & \multicolumn{3}{c}{\xrightarrow{\;\quad\qquad
                \;\,  \star \;\, \qquad\quad\,}} 
         &  {L_{}}^{\star\nts} = \ZZ[\tau]  &  \end{array}
  \renewcommand{\arraystretch}{1}
\end{equation}
Note that a CPS can be defined in the more general setting of
$\sigma$-compact locally compact Abelian groups (and beyond) for both
physical and internal space; see \cite{TAO,Moody-rev,Bernd} for more.

One of the simplest examples that needs this more general type of CPS
is the \emph{period doubling} substitution
$\varrho^{}_{\mathrm{pd}} = (ab,aa)$. It is of constant length and has
a coincidence (in the first position). It thus has pure point spectrum
by Dekking's criterion \cite{Dekking}. The corresponding CPS works
with $\RR$ as direct (or physical) space, as our guiding Fibonacci
example. However, the internal space now is $\ZZ^{}_2$, the 2-adic
integers; see \cite[Ex.~7.4]{TAO} or \cite{BMP00} for more, and for an
explicit diffraction formula.

A classic inflation rule over a ternary alphabet with a cubic
inflation factor is given by the \emph{Tribonacci} rule
$\varrho^{}_{\mathrm{Tri}} = (ab,ac,a)$, which also has a~twisted
version, namely $\varrho'_{\mathrm{Tri}} = (ba,ac,a)$. Here, we need
$\RR$ as direct and $\RR^2$ as internal space, where the windows are
now \emph{Rauzy fractals}, which are topologically regular sets
(meaning that each is the closure of its interior); see
\cite[Fig.~1]{BG-Rauzy}, and \cite{PyFo,SieThu09,Bernd} for general
background.

At this point, one might ask what happens if one extends the alphabet
to an infinite one. Like in the case of Markov chains, things get more
involved, but in the topological setting of \emph{compact alphabets},
some systematic answers are possible; see \cite{Neil} and references
therein for an introduction. One important insight is that such a step
produces many new phenomena, the perhaps most spectacular of which is
the occurrence of inflation factors that need not be an algebraic
integer, and can even be transcendental~\cite{FrGaMa22}. However, no
analogue of the projection approach is known for this generalisation,
and is not likely to exist in any obvious form.

\begin{figure}
\centering
\includegraphics[width=0.47\textwidth]{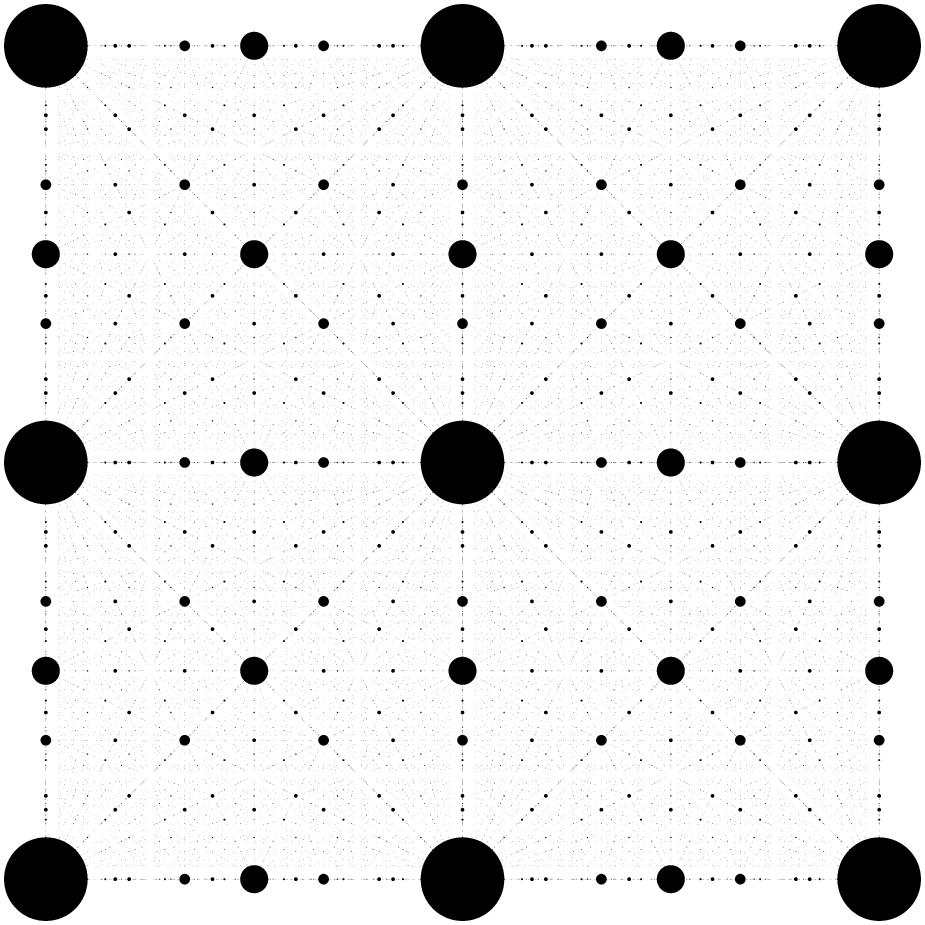}
\caption{\small Diffraction of the visible lattice points.
  They have pure point (or Bragg) diffraction. The Bragg
  peaks are represented by disks whose area is the
  height (or intensity) of the peak, located at the centre
  of the disk; see \cite[Thm.~10.5]{TAO} for more.
  \label{fig:vis} }
\end{figure}

An example of a totally different kind is provided by the square-free
integers,
\begin{equation}
  V^{}_2 \, = \, \Bigl\{ x\in\ZZ : 
      \mbox{\small \begin{tabular}{cc}$x$ is not divisible by the \\
              square of any prime\end{tabular} } \Bigr\} .
\end{equation}
It has holes of arbitrary size, and is thus \emph{not} a Delone set.
Furthermore, it cannot be turned into a Delone set by adding points
with zero density. Nevertheless, $V^{}_{2}$ still has pure point (or
Bragg) diffraction \cite{BMP00}. What is more, it can also be
described as a cut-and-project set, this time with a~compact Abelian
group $H$ as internal space. The corresponding window is a compact
subset of $H$ that has no interior, but otherwise almost everything
works as usual \cite{RS17}. A planar counterpart is the set of visible
lattice points of $\ZZ^2$, that is,
\begin{equation}
  V \, = \, \{ (m,n) \in \ZZ^2 : \gcd(m,n)=1 \} \ts .
\end{equation}
It also admits a description via a suitable CPS, and leads to the
diffraction image of Figure~\ref{fig:vis}. For details on \emph{weak
  model sets}, where one enlarges the class of admissible windows, see
\cite[Sec.~10.4]{TAO} as well as \cite{BMP00,Bernd} and references
given there.

\section{Variations and complications}\label{sec:compli}

Both $\rho^{}_1$ and $\rho^{}_2$ possess the substitution matrix $M$
from \eqref{eq:matM}, and they are the only ones compatible with
$M$. This is deceptively simple, as we can see from
\begin{equation}
  M^2 \, = \, \begin{pmatrix}
	2 & 1 \\ 1& 1  \end{pmatrix}.
\end{equation}
This matrix is compatible with precisely six substitution rules,
namely
\begin{equation}
\begin{split}
  {\rho^{}_1}^2 = (aba,ab) \ts , &
    \quad \;\: {\rho^{}_2}^2 = (aba,ba) \ts , \\
    \rho^{}_1\rho^{}_2 = (aab,ab) \ts , &
    \quad \rho^{}_2\rho^{}_1 = (baa,ba) \ts ,
\end{split}
\end{equation}
which again define the Fibonacci system, together with the two
previously mentioned rules
\begin{equation}
  (aab,ba) \quad \mbox{and} \quad (baa,ab) \ts .
\end{equation} 
They define an \emph{enantiomorphic} (or mirror) pair of different
tiling systems. Since the tilings from both systems contain the patch
$bb$, unlike the Fibonacci system, they are no longer Sturmian, and
hence more complex. Still, they admit self-similar tilings with the
same intervals as used for the Fibonacci case, and thus with the same
return module and the same CPS as above in
Eq.~\eqref{eq:CPS}. Miraculously, they are also regular model sets,
but with a much more complicated window. It is a natural question
which condition would guarantee that the window is an interval. This
has been studied extensively, and we refer readers to
\cite{BEIR07,BFS12,Cant03}.

Considering $(aab,ba)$, the \emph{reshuffled} Fibonacci substitution,
one finds a particular window pair $W_a,W_b$ of (genuine) Rauzy
fractals with a fractal boundary, here with Hausdorff dimension
\begin{equation}
  d_{\mathrm{H}} \, =\, \frac{\log (1+\sqrt{2})}{2\log (\tau)}
  \, \approx \, 0.915{\ts}785{\ts}46\dots \ts ;
\end{equation}
see Figure~\ref{fig:frac_window} below.  Its partner system has
windows $\widetilde{W}_i$ that are translates of the reflected
windows, $-W_i$, as expected. A~simple calculation reveals
\begin{equation}
  \widetilde{W}_a \,=\, \tau-1-W_a, \qquad
  \widetilde{W}_b \,=\, -1-W_b  \ts .
\end{equation}

One can now imagine how this kind of complication might grow with the
size of the matrix elements of the substitution matrix and even more
so with the size of the alphabet. This is one of the reasons why the
\emph{Pisot substitution conjecture} for alphabets with more than two
letters is still open \cite{ABBLS15}. Another reason is the
unavoidability of Rauzy fractals even for the simplest inflations once
the PF eigenvalue is an algebraic integer of degree three or higher, as
shown in \cite[Prop.~2.35]{Pleas00}.

\begin{figure}[ht]
 \includegraphics[width=0.45\textwidth]{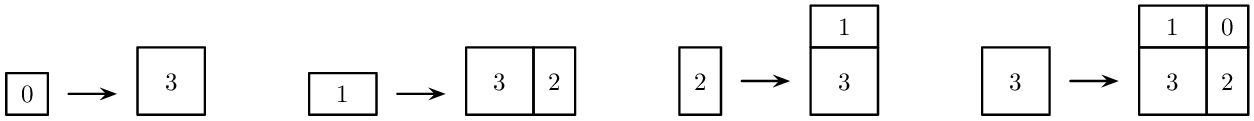}
 \caption{The Fibonacci direct product inflation rule in the
   Euclidean plane. Using the lower left corners as control points,
   one can extend the projection description of the Fibonacci chain
   to this case by doubling all dimensions and taking direct products
   of the 1D windows as the new windows here; see \cite{BFG} for details.
   \label{fig:DPV}}
\end{figure}

Let us briefly explore what else can emerge as soon as we look for
tilings of the plane. The direct product of two Fibonacci inflations
can be encoded as shown in Figure~\ref{fig:DPV}, and nothing
unexpected happens; see \cite{Clement,Laby,Ron} for some aspects and
applications.
But this is only the simplest of altogether $48$ possibilities to
define a self-similar inflation with these tile shapes (and this
inflation factor), but variations in the internal decomposition. As it
turns out \cite{BFG,BGM-Fib}, all of them are regular model sets,
though some have windows of Rauzy fractal type, by which we mean that
they have a fractal boundary. This adds another layer of complications
one has to deal with, and the understanding of the possible windows
with fractal boundary is far from complete. This kind of analysis
opens the study of higher-dimensional inflation tilings. Clearly,
there are also several variants of the famous Pisot conjecture
\cite{ABBLS15}, but a better understanding of the geometric
constraints is required for future progress.

Further variations and generalisations of the Fibonacci tiling use
another, very general approach to substitution tilings and their
relatives, which is called \emph{fusion}; see \cite{FrSa14} for a
detailed introduction. In the same paper, the \emph{scrambled
  Fibonacci tiling} was introduced. This tiling still shares a lot of
properties with the usual Fibonacci tiling, but it also serves as a
counterintuitive example of pure-point diffractive structure which is
not a Meyer set, as shown in \cite{KeSa14}. Also, multi-scale variants
have been considered, as well as concatenations of different rules
(under the name \emph{$S$-adic substitutions}); we refer to
\cite{Neil} and references therein for more.

\begin{figure*}[ht]
  \centering
  \includegraphics[width=0.95\textwidth]{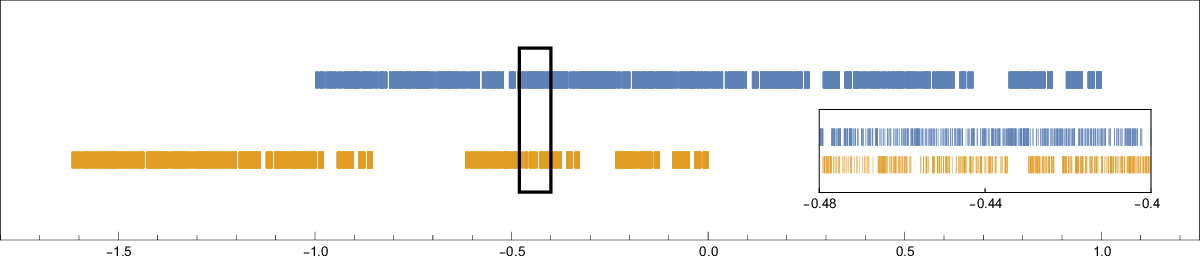}
  \caption{\small The two windows (blue/top for control points of type
    $a$ and yellow/bottom for type $b$) for the tiling given by the
    substitution $(aab,ba)$. The inlay shows a stretched view of the
    marked region. We emphasise that the sets $W_a$ and $W_b$ are
    measure-theoretically disjoint, see for example
    \cite[Cor.~6.66]{Bernd}, even though this is difficult to
    illustrate due to the high Hausdorff dimension of the boundaries.
    \label{fig:frac_window} }
\end{figure*}
 
Yet another direction was opened by \emph{random substitutions} and
\emph{inflations}. The simplest example is based on the Fibonacci
rules and can be given as
\begin{equation}
  \begin{split}
  \varrho^{}_{p} \, & = \,
  (w_{p},a) \quad \mbox{together with} \\[1mm]
  w_p \, & = \, \begin{cases}
    ab, & \mbox{with prob. } p,\\
    ba, & \mbox{with prob. } 1-p, \end{cases}
  \end{split}
\end{equation}
where the choice between $ab$ and $ba$ is randomly made at \emph{every
  step} and \emph{position}. This was introduced in \cite{GL89}, but
largely neglected for a long time. A typical realisation of this rule,
again with intervals of length $\tau$ and 1 as before, leads to
a~diffraction measure of mixed type, with pure point and absolutely
continuous contributions; see \cite{MM-diss} for further details and
illustrations.

Higher-dimensional examples are harder to find, due to geometric
constraints, but have also been analysed (some already in
\cite{GL89}), and later systematically searched for in \cite{GKM14}.
Let us add that such models are relevant, because the vast majority of
real-world quasicrystals contain a significant amount of disorder, and
are likely entropically (rather than energetically) stabilised.

\section{Equidistribution and ergodic aspects}\label{sec:ergodic}

Now, we ask what benefit we get from knowing that the Fibonacci point
set $\vL$ is a model set, $\vL = \scalebox{1.3}{$\curlywedge$}(W)$ say
with $W=(-1,\tau-1]$. The crucial observation is that $\vL^{\star}$ is
\emph{uniformly distributed} in $W$. To make sense of this statement,
one has to turn the point set $\vL$ into a natural sequence, which is
usually done by numbering its elements according to their distance
from $0$, so we write $\vL = \{x_n \, : \, n\in\NN \}$ with
$|x_n| \leqslant |x_{n+1}|$ for all $n\in \NN$. Then, the sequence
$(x^{\star}_{n})^{}_{n\in\NN}$ is uniformly distributed in $W\!$,
meaning that
\begin{equation}
  \myfrac{1}{N} \sum_{n=1}^{N} \mathbf{1}^{}_U (x^{\star}_n) \,
  \xrightarrow{\ts N \to \infty \ts } \, \vol(U)
\end{equation}
holds for every subset $U\subseteq W\nts$ such that the characteristic
function $\mathbf{1}^{}_{U}$ is Riemann integrable.  The remaining
freedom to arrange $\vL$ (when two elements have the same absolute
value) is immaterial. The corresponding property holds for $\vL_a$ and
$\vL_b$, where $\vL = \vL_a \, \dot{\cup} \, \vL_b$, relative to the
windows $W_a$ and $W_b$; see \cite{Moody} for a general account and
further references, and \cite{KN} for a systematic treatment of
equidistribution.

This property permits the determination of many frequencies coming
from averages by calculating simple integrals. In other words, we have
\emph{ergodicity}. In almost all early papers, this property was
tacitly assumed, though a~proof came much later \cite{Schl98,Moody},
and these days is a~consequence of some dynamical systems theory. The
simplest application is the determination of the relative frequencies
of points in $\vL$ of type $a$ and $b$, which gives
$\tfrac{\vol(W_a)}{\vol(W)} = \tau-1$ and
$\tfrac{\vol(W_B)}{\vol(W)} = 2-\tau$, respectively. Clearly, we know
this already from the PF right eigenvector of $M$, which seems equally
easy. However, as soon as one proceeds to the calculation of general
patch frequencies, the inflation method gets tedious, while the
uniform distribution approach often remains straightforward and easily
computable.

Let us explain this in some more detail for the Fibonacci tiling,
formulated in terms of tiles (intervals) with their left endpoints as
control points. They are all of the form
\begin{equation}
  \begin{split}
   [m+n\tau, m + (n+1)\tau] \ts , & \quad \mbox{for type}\ a, \\
   [m+n\tau, m+1 + n\tau] \ts , & \quad \mbox{for type}\ b,
  \end{split}
\end{equation}      
with $m,n \in \ZZ$, because $\vL \subset \ZZ[\tau]$ by construction
(and every other Fibonacci tiling is a translate of one with this
property). Now, each such interval has its counterpart in internal
space (as also shown in Figure~\ref{fig:fiboshear} below), sometimes
called the corresponding \emph{atomic hypersurface}, namely
\begin{align}
	\begin{split}
	m+n(1-\tau) &+ [1-\tau,\,2-\tau) \ts , \\
	m+n(1-\tau) &+ [2-\tau,\,1) \ts .
	\end{split}
\label{eq:wind_tiles}
\end{align}
Each can be understood as a coding window for the occurrence of
a specific interval in the tiling.

So, if we are given a finite set of tiles (adjacent or not), we
can decide on their joint legality in a Fibonacci point set within
$\ZZ[\tau]$, and also determine the relative patch frequency as
follows. Let $T_1,\dots \, , T_n$ be these tiles, and
$V_1,\dots \, , V_n$ their coding windows according to
\eqref{eq:wind_tiles}. Then, we consider $V_1 \cap \cdots \cap V_n$,
and obtain
\begin{equation}
  \mathrm{freq}\,(T_1,\dots \, , T_n) \, = \,
  \myfrac{\vol(V_1 \cap \cdots \cap V_n)}{\vol(W)} \ts ,
\end{equation}
which is $0$ whenever the patch is illegal (in the sense that it
cannot occur in \emph{a single} Fibonacci point set). With
$V_i = [\alpha_i,\beta_i)$, we simply get
\begin{equation}
\begin{split}
  \vol(V_1 & \cap \cdots \cap V_n) \, = \\
  & \max \{ \min_i (\beta_i ) - \max_j (\alpha_j ) \ts , 0 \} \ts ,
\end{split}  
\end{equation}
which is easy to implement. 

This approach has analogues in higher dimensions, where the inflation
method is quickly becoming impractical. In \cite{Maz}, based upon the
dualisation method from \cite{BKSZ90,KS89}, the procedure is explained
for the rhombic Penrose tiling and for the Ammann--Beenker tiling,
where exact results are derived also for several large patches.  The
patch frequencies obtained this way have interesting applications in
the theory of (discrete) Schr\"{o}dinger operators on those tilings
\cite{DEFM22}, in particular in connection with the support of
localised eigenstates.

The \emph{frequency module} of the Fibonacci tiling, which is
$\ZZ[\tau]$, is the $\ZZ$-module of rank $2$ generated by the relative
frequencies of words of length $2$ in the infinite Fibonacci word. It
is not only helpful for patch frequencies, but also appears in the
theory of one-dimensional aperiodic and ergodic Schr\"{o}dinger
operators. Indeed, consider
\begin{equation}
  \bigl(H\psi \bigr)(n) \, = \,
  \psi(n+1)+\psi(n-1)+V(n)\,\psi(n)
\end{equation}
which defines a self-adjoined operator on the Hilbert space
$\ell^2(\ZZ)$, with a potential function~$V$ that takes two values
according to the Fibonacci chain; see \cite{DGY} for a detailed survey
of the Fibonacci Hamiltonian. Then, its \emph{integrated density of
  states} (IDS) is a devil's staircase with plateaux where the IDS
takes values from the frequency module. This is a topologically rigid
structure that can be understood by Bellissard's \emph{gap labeling
  theorem}; see \cite{BIST89,BelBovGh92,BGJ93,kel} for details. Many
open questions exist around this and related topological quantum
numbers; see \cite{Johannes} for a survey.

\section{Pair correlations}\label{sec:pair}

After this fairly general description of frequencies, let us look into
the \emph{pair correlations} in more detail. For this, let
$\nu^{}_{\alpha\beta}(z)$ be the relative frequency of a tile (or
point) of type $\alpha$ and one of type $\beta$ occurring at distance
$z\in \RR$ within $\vL$. Clearly, this can only be non-zero for
$z\in \vL - \vL \subset \ZZ[\tau]$, where
\begin{equation}
  A-B \, = \, \{x-y \ : \ x\in A, \, y\in B \}
\end{equation}
is the \emph{Minkowski difference} of $A$ and $B$. In fact,
$\nu^{}_{\alpha\beta}(z)$ is positive if and only if
$z \in \vL_{\alpha} - \vL_{\beta}$, and vanishes otherwise. This is a
consequence of $\vL$ being a repetitive Delone set.

\begin{figure}[ht]
\centering
\begin{tikzpicture}[x=0.75pt,y=0.75pt,yscale=-0.8,xscale=0.8]
\draw [line width=0.75]    (95,280) -- (457,280) ;

\draw [shift={(460,280)}, rotate = 180] [fill={rgb, 255:red, 0; green,
	0; blue, 0 } ][line width=0.08] [draw opacity=0] (8.04,-3.86) --
	(0,0) -- (8.04,3.86) -- (5.34,0) -- cycle ;
	
\draw [line width=0.75]    (230,305) -- (230,83) ;

\draw [shift={(230,80)}, rotate = 90] [fill={rgb, 255:red, 0; green,
	0; blue, 0 } ][line width=0.08] [draw opacity=0] (8.04,-3.86) --
	(0,0) -- (8.04,3.86) -- (5.34,0) -- cycle ;
		
\draw [thick]    (350,280)--(230,130) -- (110,280) ;
\draw [thick,red]   (303.2,280)--(230,188.5) -- (156.8,280) ;
\draw [thick, dashed]   (423.2,280)--(350,188.5)--(303.2,188.5) -- (230,280) ;
	
\draw (137.15,196) node [anchor=north west][inner sep=0.75pt]    {$g^{}_{aa}$};
\draw (165.95,219) node [anchor=north west][inner sep=0.75pt]    {$g^{}_{bb}$};
\draw (377.15,200) node [anchor=north west][inner sep=0.75pt]    {$g^{}_{ba}$};
\draw (345.95,300) node [anchor=south west][inner sep=0.75pt]    {\tiny $1$};
\draw (420.35,300) node [anchor=south west][inner sep=0.75pt]
	{\tiny $\tau $};
\draw (293,300) node [anchor=south west][inner sep=0.75pt]
	{\tiny $\tau^{-1}$};
\draw (233,113) node [anchor=north west][inner sep=0.75pt]
	{\tiny $\tau^{-1}$};
\draw (233,175) node [anchor=north west][inner sep=0.75pt]
	{\tiny $\tau^{-2}$};	
\end{tikzpicture}
\caption{\small The continuous functions $g^{}_{\alpha \beta}$
  describing the pair correlations in internal
  space. \label{fig:functions}}
\end{figure}
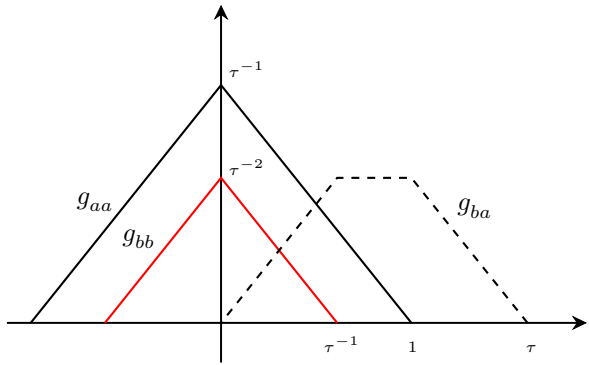

One can now use the method explained above. It can be simplified
by observing that $ z\in\vL_{\alpha} - \vL_{\beta}$ is equivalent to
$z\in \ZZ[\tau]$ together with $z^{\star} \in W_{\alpha} -
W_{\beta}$. Calculating the frequencies leads to
\begin{equation}
  \nu^{}_{\alpha\beta}(z) \, = \, g^{}_{\alpha \beta}(z^{\star})
\end{equation}
for $z\in\ZZ[\tau]$, with $g^{}_{\alpha \beta}$ the simple continuous
functions shown in Figure~\ref{fig:functions}.  Explicitly, they read
\begin{align}
  g^{}_{aa}(y) & = \max \bigl( \tfrac{1-|y|}{\tau}, 0 \bigr),
   \nonumber \\
  g^{}_{ba}(y) & = \begin{cases} \frac{y}{\tau},
    & \mbox{if } 0\leqslant
    y \leqslant \tau-1, \\
    \frac{1}{\tau^2}, & \mbox{if } \tau-1\leqslant y \leqslant 1, \\
    1-\frac{y}{\tau}, & \mbox{if } 1\leqslant y  \leqslant \tau, \\
    0, & \mbox{otherwise,} \end{cases} \\
  g^{}_{bb}(y) & = \max \bigl(\tfrac{1-\tau \ts |y|}{\tau^2},
            0 \bigr), \nonumber
\end{align}
together with $g^{}_{ab} (y) = g^{}_{ba} (-y)$.

The total pair correlation (also known as the \emph{autocorrelation}
of $\vL$) is
\begin{equation}\label{eq:def_autocor}
  \nu(z) \, = \sum_{\alpha, \beta} \nu^{}_{\alpha\beta}(z)
  \, = \, g(z^{\star}).
\end{equation}
Here,
$g(y) = \frac{1}{\tau} \bigl( \mathbf{1}^{}_{-W} \ast \mathbf{1}^{}_{W}
\bigr) (y)$ is the \emph{covariogram} of the total window, and similar
representations hold for the $\nu^{}_{\alpha\beta}$, namely
\begin{equation}
  \nu^{}_{\alpha\beta}(z) \, = \, \myfrac{1}{\tau}
  \bigl( \mathbf{1}^{}_{-W_{\alpha}}
  \ast  \mathbf{1}^{}_{W_{\beta}} \bigr) (z^{\star}) \ts ,
\end{equation}
which explains the result shown in Figure~\ref{fig:functions}. 

How is all this reflected in the inflation picture? As recently shown
in \cite{BG,BGM}, the correlation coefficients satisfy the exact
\emph{renormalisation relations}
\begin{equation}\label{eq:renorm}
\begin{split}
  \nu^{}_{aa}(z)&=\tfrac{1}{\tau}\big(\nu^{}_{aa}
  \big(\tfrac{z}{\tau}\big) +\nu^{}_{ab}\big(\tfrac{z}{\tau}\big) \\
  & \qquad \quad +\nu^{}_{ba}\big(\tfrac{z}{\tau}\big) +
    \nu^{}_{bb}\big(\tfrac{z}{\tau}\big)  \big), \\
    \nu^{}_{ab}(z)&=\tfrac{1}{\tau}\big(\nu^{}_{aa}
    \big(\tfrac{z}{\tau}-1\big)  +\nu^{}_{ba}
    \big(\tfrac{z}{\tau}-1\big)   \big), \\
    \nu^{}_{ba}(z)&=\tfrac{1}{\tau}\big(\nu^{}_{aa}
    \big(\tfrac{z}{\tau}+1\big)  +
    \nu^{}_{ab}\big(\tfrac{z}{\tau}+1\big)   \big), \\
    \nu^{}_{bb}(z)&=\tfrac{1}{\tau}\nu^{}_{aa}
    \big(\tfrac{z}{\tau} \big),
\end{split}
\end{equation}
with $z\in \vL-\vL$. 

This is an infinite set of linear equations. A~\emph{finite} subset of
them closes, namely the ones with $|z| \leqslant \tau$ on the left
side. Subject to the constraints on the possible $z\in \vL - \vL$,
this subset has a one-dimensional solution space, while all remaining
coefficients are recursively determined from the ones of this
subset. Specifying $\nu^{}_{aa}(0) + \nu^{}_{bb}(0) = 1$ then gives
the solution described above in terms of the $g$-functions.

No similarly simple renormalisation seems to exist for $\nu(z)$. But
one can use \eqref{eq:def_autocor} together with \eqref{eq:renorm}
iteratively to derive the relation
\begin{equation}
  \begin{split}
  \nu(z) \, &= \, \myfrac{1}{\tau^2} \, \nu \Bigl(\myfrac{1}
              {\tau^2} \Bigr) \, + \\[1mm]
    &\sum_{n\in\ZZ} \myfrac{1}{\tau^{|n|+1}} \, \nu \biggl(
      \myfrac{z + \mathrm{sgn}(n) \bigl(({-}\tau)^{|n|} -1 \bigr)}
      {\tau^{|n|+1}} \biggr) 
  \end{split}
\end{equation}
which can be interpreted in terms of the functions
$g^{}_{\alpha\beta}$ and rescaled/translated versions of them.
The treatment via the renormalisation relations provides a powerful
tool for the cases where the covariogram is difficult to access, for
instance when the windows have fractal boundaries.
 
Now that we know the correlation coefficients for the self-similar
Fibonacci tiling, it is an obvious question whether (and how) one can
also get them for modified versions, in particular for the case that
we use intervals of two arbitrary lengths. This is possible as long as
the average interval length is $\tfrac{\sqrt{5}}{\tau}$, which is the
one from our self-similar case. Other situations can be obtained from
here by a simple global rescaling; see \cite{Jan} for more. The same
approach also works for primitive inflation tilings in higher
dimensions. The exact renormalisation relations have been used to rule
out absolutely continuous spectral contributions to the diffraction in
various examples, including the Godr\`{e}che--Lan\c{c}on--Billard tiling
\cite{BGM}; see \cite{TAO} for background and further
references. Indeed, its diffraction is purely singular continuous
(except for the trivial Bragg peak at $k=0$), which had been the common
assumption for this non-PV inflation rule.

\section{Shape changes}\label{sec:shape}

\begin{figure}[ht]
  \includegraphics[width=0.45\textwidth]{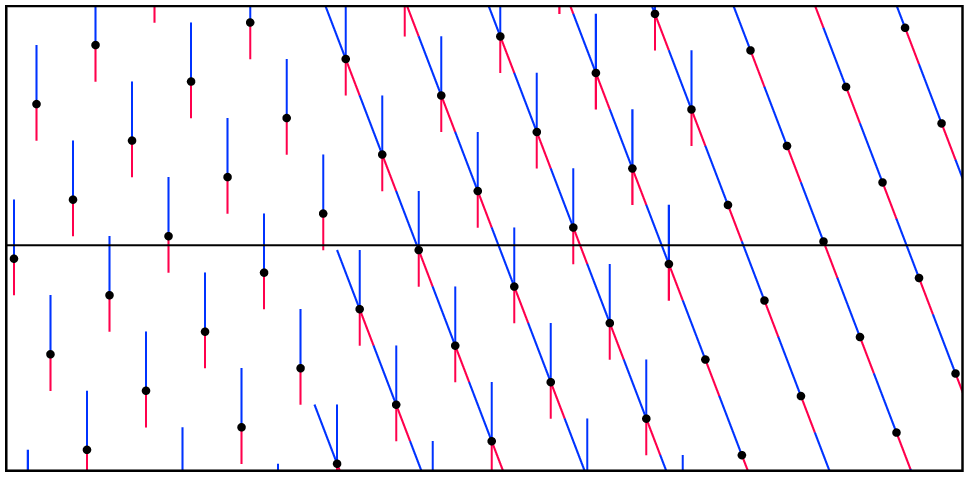}
  \caption{\small Original and sheared Fibonacci tiling, obtained as a
    section through a periodic array of \textit{atomic
      hypersurfaces}. \label{fig:fiboshear}}
\end{figure}

To illustrate the effect of changing the relative tile lengths of the
Fibonacci tiling, it is most convenient to view it as a
\textit{section} through a periodic array of atomic hypersurfaces (see
Figure~\ref{fig:fiboshear}). This is equivalent to the cut-and-project
construction, as explained in detail in \cite[Sec.~7.5.1]{TAO}. If we
shear these atomic hypersurfaces \textit{parallel} to the cut
direction, maintaining the lattice sites at which they are attached,
no vertices appear or disappear, nor change their type. Also, the
periodicity of the array of hypersurfaces remains the same. The only
effect is a~change in the relative tile lengths. In
Figure~\ref{fig:fiboshear}, we illustrate the change from original
tile lengths (left) to all tile lengths equal (right). In the middle,
both the sheared and unsheared hypersurfaces are displayed.  As one
can see, the deformed tiling is obtained by projecting the same
lattice sites, which fall into the same window, along a new projection
direction parallel to the sheared atomic hypersurfaces.

The amount by which a vertex is moved depends non-locally on its
environment. This length change is thus not locally derivable from the
undeformed tiling, but it still induces a (non-local) topological
conjugacy of the dynamical system, because it preserves the
higher-dimensional lattice, and thus commutes with the translation
action. In other words: the deformation does not mess up the aperiodic
translational order. In fact, one can show \cite{CS1,CS2} that, as
long as the overall scale is maintained, the Fibonacci tiling does not
admit length changes which affect the dynamical system in a relevant
way, which is a topological property of the Fibonacci hull; see
\cite{AP,tilingsbook} for background on the topological methods applied
here.

In the Hat and the Spectre tiling, the very same mechanism is at work.
Both are obtained as shape changes from a cut-and-project tiling (with
suitable control points), with copies of $\RR^2$ as physical and
internal spaces, and a four-dimensional lattice \cite{BGS}. This
explains their pure point diffraction property (which was numerically
calculated in \cite{Soc}), and their role as mathematical
quasicrystals. The interesting feature in comparison to previous
aperiodic monotiles is that they are quasiperiodic rather than
limit-periodic (both in the sense of \emph{mean almost periodicity},
compare \cite{LSS}, which is a recent extension of the notion of
almost periodic functions in a measure-theoretic setting).  It will be
interesting to see which other tilings of this kind will be
discovered, in particular in more than two dimensions.

\section{Diffraction}\label{sec:diffract}

Let us put all $\nu(z)$ in one object, 
\begin{equation}
   \gamma \, \defeq \sum_{z\in \vL - \vL} \nu(z) \, \delta_{z} \ts ,
\end{equation}
where $\delta_{z}$ is the \emph{Dirac measure} (or distribution) at
$z$. It is defined by $\delta_z(\phi) = \phi(z)$ for any function that
is continuous at $z$. Here, $\gamma$ is a \emph{measure}, and should
be considered as the (rigorous) infinite-size analogue of the
Patterson function of crystallography. It is the natural
autocorrelation measure of
\begin{equation}
  \omega \, = \, \delta^{}_{\!\vL} \, \defeq
  \sum_{x\in \vL} \delta_x \ts ,
\end{equation}
the Dirac comb of the point set $\vL$. The autocorrelation measure is
usually defined as
\begin{equation}
  \gamma \, =  \lim_{n\to \infty} \myfrac{1}{2n}
  \sum_{x,y\in\vL^{(n)}} \delta^{}_{x-y} \ts ,
\end{equation}
with $\vL^{(n)} = \vL \cap [-n,n]$. The existence of the limit is a
consequence of the underlying ergodic properties of $\vL$.  Here,
$\gamma$ is a strongly almost periodic measure with Fourier transform

\begin{equation} 
   \widehat{\gamma} \, = \sum_{k \in L^{\circledast}} I(k)\, \delta^{}_k
\end{equation}
which is the \emph{diffraction measure} of $\omega$. The supporting
set is $L^{\circledast} = \ZZ[\tau]/\sqrt{5}$, also known as the
dynamical or the Fourier--Bohr spectrum; see
Section~\ref{sec:dynamical} for further details. In crystallography,
this is often called the \emph{reciprocal lattice}, though it is not a
lattice in the mathematical sense once the structure is
non-periodic. The spectrum satisfies $L^{\circledast} = \pi (\cL^*)$,
where $\cL^*$ is the dual of the embedding lattice $\cL$. The
intensities at $k\in L^{\circledast}$ are given by $I(k) = |A(k)|^2$
with the amplitudes (or \emph{Fourier--Bohr coefficients})
\begin{align}
  A(k) \, &= \lim_{n\to \infty}\myfrac{1}{2n}
            \sum_{x\in\vL^{(n)}} \ee^{-2\pi\ii kx}
    \nonumber \\[1mm]
          &= \, \begin{cases}
            \myfrac{\dens(\vL)}{\vol(W)} \;
            \widehat{\mathbf{1}^{}_W}(-k^{\star}) , &
            \mbox{if} \ k\in L^{\circledast}, \\
            0, \ & \mbox{otherwise}, \end{cases}
    \label{eq:FB}
\end{align}
where the limit always exists; see \cite{BH} for an elementary proof
of this formula. Note that the density of $\vL$ is well defined, and
that $\widehat{\ts\mathbf{1}^{}_{W}\ts}$ is the Fourier transform of
the characteristic function of $W$.

The importance of the FB coefficients can hardly be overstated. They
are also instrumental in the recent classification of pure point
diffraction via almost periodicity \cite{LSS,LSS23}, and they play an
important role in the dynamical systems approach to aperiodic order,
as we shall see in Section~\ref{sec:dynamical}.

More generally, one is interested in weighted Fibonacci combs, such as
\begin{equation}
  \omega \, = \, h^{}_a \ts \delta^{}_{\!\vL_a} +
     h^{}_b \ts \delta^{}_{\!\vL_b} \ts ,
\end{equation}
which gives scattering strength $h_\alpha$ to points of type
$\alpha$. Now, the autocorrelation becomes
\begin{equation}
  \begin{split}
  \gamma \, & \, =  |h^{}_a|^2 \, \gamma^{}_{aa} +
              h^{}_a\,\overline{h^{}_b}\, \gamma^{}_{ab} \\[1mm]
            & \qquad + \overline{h^{}_a} \, h^{}_b \,
              \gamma^{}_{ba} + |h^{}_b|^2 \, \gamma^{}_{bb}
  \end{split}
\end{equation}
with diffraction 
\begin{equation} 
  \widehat{\gamma} \, =  \sum_{k \in L^{\circledast}}
  \bigl| h^{}_a A^{}_a (k) + h^{}_b A^{}_b(k) \bigr|^2 \,
  \delta^{}_k \ts ,
\end{equation}
where $A_{\alpha}(k)$ is the FB coefficient of
$\delta^{}_{\!\vL_{\alpha}}$ at $k$. The formula reflects the phase
consistency property as proved in \cite{BGM,LSS}. More general
weighting schemes can be considered, as outlined in \cite{Nicu}.
  
\begin{figure}[ht]
\centering
\includegraphics[width=0.45\textwidth]{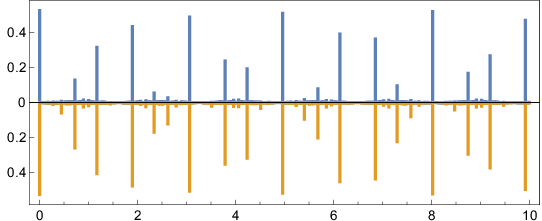}
\caption{\small The diffraction image for the Fibonacci tiling
  (yellow, bottom) and for the reshuffled Fibonacci (blue, top) for
  $0\leqslant k \leqslant 10$. The intensity of every Bragg peak (the
  absolute square of the Fourier--Bohr coefficient) is given by the
  height of the bar. The intensity of the central peak at $0$ is
  $(\tau+1)/5 \approx 0.5236$. }
\label{fig:diffraction}
\end{figure}

In Figure~\ref{fig:diffraction}, we show the diffraction of the
Fibonacci chain in comparison to the one from the reshuffled
version. The latter has the much more complicated windows from
Figure~\ref{fig:frac_window} with their fractal boundaries. In fact,
calculating the FB coefficients for them requires a method from
\cite{BG-Rauzy} to compute $\widehat{\mathbf{1}^{}_W}$, which is based
on the \emph{Fourier matrix cocycle}. Note that a numerical approach
via the Fourier transform of large finite patches converges in
principle, but rather slowly.

This is actually a typical situation, which is not restricted to
windows in one dimension.  It also arises in the case of the (twisted)
Tribonacci inflation $\varrho^{(\prime)}_{\mathrm{Tri}}$ mentioned
earlier, where the windows are two-dimensional Rauzy fractals. In the
standard case, the windows are simply connected, and the numerical
approach still works reasonably well, while it fails rather badly for
the twisted case, where the windows are `spongy'. This was studied in
more details for the plastic number inflation in \cite{BG-plastic},
which has windows of a similar type. Generally, the difficulty of a
reliable (numerical) calculation increases with the Hausdorff
dimension of the boundary.

\begin{figure}[t]
\centering
\includegraphics[width=0.45\textwidth]{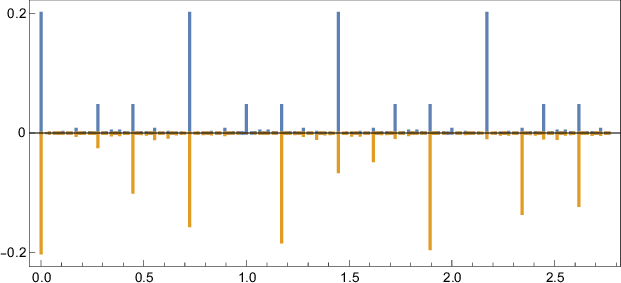}
\caption{\small Diffraction of the symbolic Fibonacci chain with equal
  tile lengths of $\sqrt{5}/ \tau$ and weights $h_a=1$, $h_b=0$ (blue,
  top) and in comparison to its geometric realization with natural
  tile lengths $\tau$ and $1$ (yellow, bottom). The intensity of the
  central peak at $0$ is $1/5$ and the peaks are shown for wave numbers
  in $[0, 2\sqrt{5}/ \tau)$, which corresponds to two fundamental cells
  of the average lattice underneath the symbolic case.
  \label{fig:diff_symb_vs_geom} }
\end{figure}

Let us mention that the embedding formalism also permits to compute
the diffraction of the Fibonacci chain under the shape changes
illustrated in Figure~\ref{fig:fiboshear}. This results in
\emph{deformed} model sets, compare \cite[Ex.~9.9]{TAO} for a related
example and \cite{BerDun00} for a general discussion, which still
leads to a closed formula for the diffraction.
Figure~\ref{fig:diff_symb_vs_geom} illustrates the result for the
extreme case that the interval lengths become equal. The aperiodicity
is still present (via two different weights for the point types) and
clearly visible from the peaks, despite the fact that the (coloured)
points now live on a lattice. In this case, in line with
\cite[Thm.~10.3]{TAO}, the diffraction measure is periodic.  The
analogous situation occurs for the Hat tiling \cite{BGS}, then with a
hexagonal lattice. The details of this structure are complicated by
the fractal nature of the windows, and will be analysed more
extensively in \cite{BMM}.

\begin{figure*}
\centering
\includegraphics[width=0.45\textwidth]{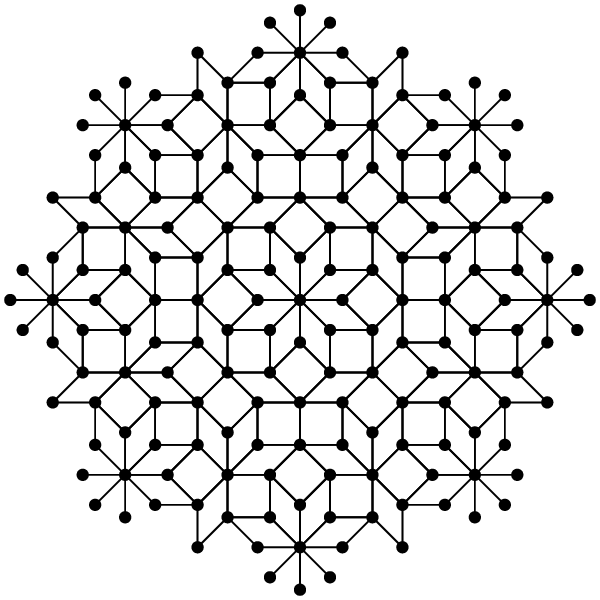}\hspace*{2em}
\includegraphics[width=0.48\textwidth]{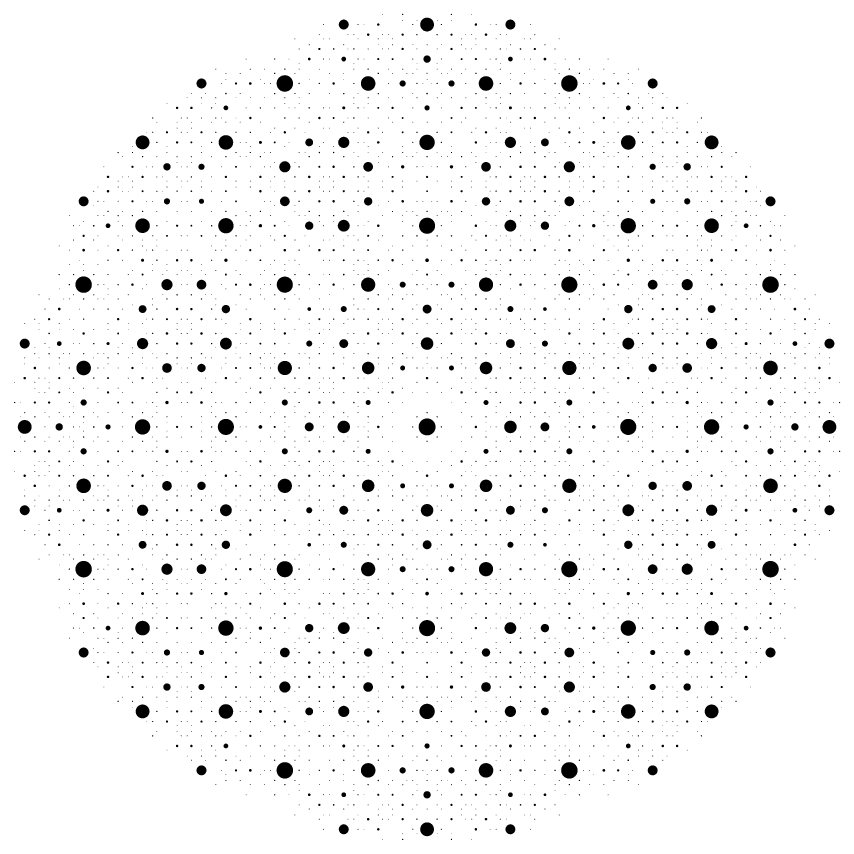}
\caption{\small The left panel shows a finite, eightfold symmetric
  patch of the Ammann--Beenker tiling. The vertex points can be
  obtained as a model set via the hypercubic lattice and a regular
  octagon as window in internal space. The right panel displays
  the diffraction of the (infinite) Ammann--Beenker
  model set. Bragg peaks are again represented as disks, as in
  Figure~\ref{fig:vis}. \label{fig:AB} }
\end{figure*}

In two dimensions, among the most prominent examples are the rhombic
Penrose, the Ammann--Beenker (AB), and the square-triangle tiling due
to Schlottmann \cite[Sec.~6.3.1]{TAO}, with $10$-, $8$- and $12$-fold
symmetry, respectively. The CPS is simplest for the AB tiling,
producing for instance the symmetric patch and the diffraction image
of Figure~\ref{fig:AB}. This is based on the lattice $\ZZ^4$ in
$\RR^4$; see \cite[Ex.~7.8]{TAO} for details. The vertices of the
rhombic Penrose tiling fall into four distinct translation classes, as
analysed in detail in \cite{BKSZ90}, thus completing the pioneering
work by de Bruijn \cite{dBr81}. The square-triangle tiling is the most
difficult of these three, because it has a window with twelvefold
symmetry but fractal boundary. This is unavoidable for square-triangle
tilings, and one inflation tiling of this kind, due to Schlottmann, is
described in detail in \cite[Sec.~6.3.1 and Fig.~7.10]{TAO}.

\section{Dynamical systems approach}\label{sec:dynamical}

The dynamical systems mentioned earlier have always been important
objects of mathematical research, but their importance was hardly
recognised in physics. This changed in the context of aperiodic order,
where they bring important insight, in particular for an adequate
definition of symmetries. An individual tiling need not be mirror or
rotation symmetric, but its hull possesses this symmetry, which is the
correct way to analyse this. Also, many other properties do \emph{not}
depend on the individual element of the hull, but hold for the hull as
a whole, such as diffraction or other spectral features.

Let us discuss this a bit more, in the context of the Fibonacci tiling
dynamical system $(\YY,\RR)$ from Section~\ref{sec:scene}. Since $\YY$
is compact, there exists (by general arguments) at least one
probability measure that is invariant under the translation action. In
fact, in this example, there is precisely \emph{one} such measure,
$\mu_{_\mathrm{F}}$, and this turns the topological dynamical system
$(\YY,\RR)$ into a measure-theoretical one, denoted by
$(\YY,\RR,\mu_{_\mathrm{F}})$.

The measure $\mu_{_\mathrm{F}}$ is \emph{ergodic} (meaning that every
invariant subset of $\YY$ has either measure 0 or 1), and all
averages along Fibonacci chains can be written as integrals over $\YY$
with respect to $\mu_{_\mathrm{F}}$. Another connection concerns the
spectral theory of the Fibonacci chain. Given $\YY$ and
$\mu_{_\mathrm{F}}$, one can define the Hilbert space
$L^2(\YY,\mu_{_\mathrm{F}})$ of square-integrable functions
$f:\YY \rightarrow \CC$, where the inner product is 
\begin{equation}
  \langle g \ts\ts | \ts f \rangle \, \defeq
  \int_{\YY}\, \overline{g(Y)} f(Y) \dd \mu_{_\mathrm{F}}(Y) \ts .
\end{equation}

The crucial observation, due to Koopman, now is that the translation
action on the tiling space induces a family of \emph{unitary
  operators} $T^{}_{t}$ on $L^2(\YY,\mu_{_\mathrm{F}})$ via
$(T^{}_{t} f)(Y)\defeq f(Y\! {-}\ts t)$, where $Y\in \YY$ is
a~tiling and $Y\! {-}\ts t$ its translate. Indeed, for arbitrary
$f,g \in L^2(\YY,\mu_{_\mathrm{F}})$, one gets
\begin{equation}
\begin{split}
  \langle T^{}_{t} g \ts\ts | \ts T^{}_{t} f \rangle \, & = 
  \int_{\YY}\, \overline{g(Y\! {-}\ts t)} f(Y\! {-}\ts t)
  \dd \mu_{_\mathrm{F}}(Y) \\[1mm]
  & =   \int_{\YY}\, \overline{g(Y)} f(Y) \dd \mu_{_\mathrm{F}}(Y)
  \, = \, \langle g \ts\ts | \ts f \rangle,
\end{split}
\end{equation}
where the second step follows by a change of variable transform and
the translation invariance of $\mu_{_\mathrm{F}}$.

Since the $T^{}_{t}$ commute with one another for all $t\in\RR$, the
operators possess simultaneous eigenfunctions (if any), and the
remarkable property here is that there is a (countable) set of
eigenfunctions which span $L^2(\YY,\mu_{_\mathrm{F}})$. One then
says that $(\YY,\RR,\mu_{_\mathrm{F}})$ has \emph{pure point
  dynamical spectrum}. The theory of such systems was developed by
Halmos and von Neumann \cite{HvN44}, and is a cornerstone of dynamical
systems theory. Here, the connection can be made more concrete, which
gives a link to diffraction theory; see \cite{Sol97,Q,BL} and
references therein.
 
Let us select a Fibonacci Delone set $\vL$, the point set of a tiling
$Y\in \YY$, and consider the corresponding \emph{Fourier--Bohr
  coefficient} $A^{}_{\!\vL} (k) = A (k)$ as defined in
Eq.~\eqref{eq:FB} for arbitrary $k\in \RR$, often called the wave
number. In physics, this is the (complex) amplitude of the
structure. The strong ergodicity of $\mu_{_\mathrm{F}}$ implies that
the limit always exists, and one obtains the formula given in
\eqref{eq:FB}, with $W = (-1,\tau -1]$ if $\vL$ is our Fibonacci point
set from above, and $^{\star}$ is the star map. In diffraction, as
explained earlier, we get the intensity of the Bragg peak at $k$ as
$I(k) = |A^{}_{\!\vL}(k)|^2$.
  
The connection to dynamics now comes from the observation of how
$A^{}_{\!\vL}(k)$ behaves under translations of $\vL$, where we get
\begin{equation}
  A^{}_{t+\vL}(k) \, = \, \ee^{-2\pi \ii kt} A^{}_{\!\vL}(k) \ts ,
\end{equation}
for all $k\in \RR$. When $A^{}_{\!\vL}(k) \neq 0$, which is true
unless we hit an exceptional extinction point, this can be considered
as an \emph{eigenfunction equation} with eigenvalue
$\ee^{-2\pi \ii kt}$, because the left-hand side is the translate of
$A^{}_{\!\vL}$ by $t$, evaluated at the point $k$ in reciprocal space.
To avoid the appearance of extinctions, one can consider a Dirac comb
with different weights for points of type $a$ and $b$, and define the
FB coefficients accordingly. Then, generically, there are no
extinctions, and one obtains a complete set of eigenfunctions.  In our
guiding example, they turn out to be \emph{continuous} on $\YY$, and
are thus called \emph{topological eigenvalues}; see \cite{Johannes}
for interesting connections with topological invariants.

Since $\RR$ is a continuous group, it is advantageous to use
the wave numbers $k$ to label the eigenvalues. Thus,
$L^{\circledast} = \ZZ[\tau]/\sqrt{5}$ is called the \emph{dynamical
  spectrum} of $(\YY,\RR,\mu_{_\mathrm{F}})$ in additive notation.

A fundamental insight (based on Dworkin's argument \cite{Dwo93,DM08})
now is that a \emph{Delone dynamical system} (such as our Fibonacci model
set) has pure point dynamical spectrum if and only if the diffraction
measure of a typical element of the tiling hull has pure point (or
Bragg) diffraction \cite{BL,LMS,LS07}. In good cases (such as our guiding
example), every element is typical, in others (such as the visible
lattice points), one has to make the correct choice. This connection is
also instrumental in the recent analysis of pure point diffraction via
averaged versions of almost periodicity; see \cite{LSS} and references
therein for more.

\section*{Acknowledgements}
	
It is our pleasure to thank Neil Ma\~{n}ibo, Andrew Mitchell and
Lorenzo Sadun for discussions and useful hints on the manuscript.  We
are grateful to Ron Lifshitz and two anonymous referees for their
thoughtful comments, which helped us to improve the presentation.

This work was supported by the German Research Council (Deutsche
Forschungsgemeinschaft, DFG) under contract SFB-1283/2 (2021 --
317210226).  
	
\end{document}